\documentclass[a4paper,10pt]{article}

\usepackage{amsmath, amssymb, amsthm}
\usepackage[english,francais]{babel}
\newtheorem{thm}{Theorem}[section]
\newtheorem{cor}[thm]{Corollary}
\newtheorem{lem}[thm]{Lemma}

\theoremstyle{definition}

\newcommand{\R}{\mathbb{R}}

\newcommand{\ud}{\mathrm{d}}

\newcommand{\be}{\begin{equation}}
\newcommand{\ee}{\end{equation}}
\newcommand{\bs}{\begin{split}}
\newcommand{\es}{\end{split}}
\newcommand{\bee}{\begin{equation*}}
\newcommand{\eee}{\end{equation*}}

\title{Existence and conservation laws for the Boltzmann-Fermi equation in a general domain}
\author{Thibaut Allemand}

\begin{document}

\maketitle

\begin{abstract}
\selectlanguage{english}
We prove an existence theorem for the Boltzmann-Fermi-Dirac equation for integrable collision kernels in possibly bounded domains with specular reflection at the boundaries, using the characteristic lines of the free transport. We then obtain that the solution satisfies the local conservations of mass, momentum and kinetic energy thanks to a dispersion technique.

\end{abstract}

\section{Introduction}
\label{}



In this note, we are concerned with the Boltzmann equation for a gas of fermions. It has been studied for example by J. Dolbeault \cite{dolbeault}, P.-L. Lions \cite{lions} and X. Lu \cite{lu}. Our goal is double: first, proving that the existence result and the properties of the solution proved in \cite{dolbeault} can be extended to more general domains than $\R^3$, if we supplement the equation with a specular reflection condition at the boundaries; and then, establishing the local conservation laws. These two properties are needed in the study of hydrodynamic limits, especially in the incompressible inviscid regime, which is studied in a forthcoming paper \cite{moi}.

Let us detail the problem. Let $\Omega$ be a subset of $\R^3$ having bounded curvature, and regular enough such that the function $x\mapsto n(x)$, where $n(x)$ is the outer unit normal at point $x$ of $\partial \Omega$, can be extended into a continously differentiable function on $\R^3$.

The Boltzmann equation for a gas of fermions reads

\be
\label{BFD}
 \partial _t f +v.\nabla_x f = Q(f)
\ee
where $f(t,x,v)$ is the density of particles which at time $t\in\R_+$ are at point $x\in\Omega$ with velocity $v\in\R^3$. It is supplemented with an initial condition
\be
\label{initial}
 f(0,x,v)=f_0(x,v) \qquad\forall (x,v)\in \Omega\times \R^3
\ee
and a boundary condition
\be
\label{boundary}
 f(t,x,v)=f(t,x,R_x(v))\qquad\forall (x,v)\in \partial\Omega\times\R^3 ~~\textrm{such that}~~n(x).v<0
\ee
where $R_x(v)$ is the specular reflection law
\[
 R_x(v)=v-2(v.n(x))n(x).
\]

 The collision integral $Q(f)$ is given by

\be
\label{collision integral}
 Q(f)=\int_{S^2}\int_{\R^3}b(v-v_*,\omega) \left(f'f_*'(1-f)(1-f_*)-ff_*(1-f')(1-f_*') \right)\ud v_* \ud\omega
\ee
with the usual notations
\[
 f_*=f(t,x,v_*),\quad f'=f(t,x,v'),\quad f_*'=f(t,x,v_*')
\]
and where the precollisional velocities $(v',v_*')$ are deduced from the postcollisional ones by
\[
  v'=v-(v-v_*).\omega \omega,\qquad v_*'=v_*+(v-v_*).\omega \omega.
\]
The collision integral differs from the classical one by the terms $(1-f)$, which take into account the Pauli exclusion principle.

The function $b(w,\omega)$, known as the collision kernel, is measurable, a.e. positive, and is assumed to be in $L^1$ as in \cite{dolbeault}. For the study of the conservations, we will moreover assume that 
\be
\label{hyp}
b(w,\omega)=q(|w|,|w.\omega|),
\ee
which is physically relevant, and endows (formally) $Q(f)$ with symmetry properties:
\be
\label{symmetry}
 \int_{\R^3}(1,v,|v|^2)Q(f)\ud v=0
\ee
for all $f$ such that the integrals make sense. Then, integrating equation \eqref{BFD} against $1,v,|v|^2,$ we obtain that the solution $f$ satisfies formally the conservation of the local mass, momentum and kinetic energy:
\begin{equation}
\label{conslaw}
\partial_t \int_{\R^3}\begin{pmatrix} 1 \\v\\|v|^2\end{pmatrix}f\ud v + \nabla_x.\int_{\R^3}\begin{pmatrix}v\\v\otimes v\\v|v|^2\end{pmatrix}f\ud v=0.
\end{equation}

This note aims at showing first an existence theorem for the initial-boundary value (IBV) problem \eqref{BFD}-\eqref{boundary}, and then that the local conservation laws \eqref{conslaw} are satisfied rigourously by the solution of the problem.

\section{The IBV problem}

Following \cite{kaniel} and \cite{hamdache}, we define the characteristic lines of the free transport equation
\[
 \partial_t f +v.\nabla_x f=0
\]
 in the following way: for $(x,v)\in\Omega\times\R^3$, the characteristic line is given by $x+vt$, at least if $t$ is small. Then, let $t_0(x,v)$ be the first value of $t$ for which $x+t_0v \in \partial \Omega$. Then, for $t>t_0$, the trajectory continues as $x+t_0v+(t-t_0)R_{x+t_0v}(v)$. If it intersects one more time the boundary of $\Omega$ at time $t_1(x,v)$, then the trajectory continues with velocity $R_{x+t_0v+(t_1-t_0)R_{x+t_0v}(v)}(R_{x+t_0v}(v))$, and so on. We then define a family of maps $\{\Psi^t\}$ called the trajectory maps, $\Psi^t(x_0,v_0)$ beeing the point in phase space at wich we arrive at time $t$ following the trajectory line issuing from $(x_0,v_0)$. For $t=t_0,t_1...$, we define $\Psi^t$ to be continous (in time) from the right. Since $\Omega$ has bounded curvature, the trajectory intersects the boundary finitely many times in finite times. Moreover, the energy is conserved at each reflection, $\Psi^t$ maps $\bar \Omega \times \R^3$ onto itself for every $t$, and the jacobian of $\Psi^t$ is always unity \cite{kaniel}.

With this construction, we now define
\[
 f^\sharp (t,x,v)=f(t,\Psi^t(x,v)),
\]
that is, we conjugate $f$ with the free transport semigroup. This notation allows to reformulate the boundary condition in a simpler way. Indeed, notice that if $\Psi^{t-0}(x_0,v_0)=(x,v)$, then $\Psi^t(x_0,v_0)=\Psi^{t+0}(x_0,v_0)=(x,R(v))$, so that the boundary condition 
can be written
\be
\label{BC}
 f^\sharp(t-0,x_0,v_0)=f^\sharp(t+0,x_0,v_0),
\ee
that is, $t\mapsto f^\sharp(t,x,v)$ is continuous.

The existence result for the IBV problem \eqref{BFD}-\eqref{initial}-\eqref{BC} is a consequence of the following lemma:

\begin{lem}
 Let $f,h\in L^1_{\textrm{loc}}(\R_+\times\Omega\times\R^3)$. Then $f$ is a solution of
\be
\label{transport} 
\partial_t f +v.\nabla_x f =h \qquad\qquad\quad\textrm{in}~~\mathcal D'(\R_+\times\Omega\times\R^3)
\ee
 with boundary condition \eqref{BC} if and only if for almost all $(x,v)\in \Omega\times\R^3$, $f^\sharp$ is absolutely continuous with respect to $t$, $h^\sharp(t,x,v)\in L^1_{\textrm{loc}}(\R_+)$ and
\[
 f^\sharp(t,x,v)=f_0^\sharp(t,x,v)+\int_0^th^\sharp(s,x,v)\ud s.
\]

\end{lem}
\textbf{Proof.}
 It is very similar to what is done in \cite{diperna, hamdache}. It consists un multiplying \eqref{transport} by $\chi(\Psi^t(x,v)))\phi(t)$ and integrating, with $\chi \in \mathcal D(\Omega\times \R^3)$ and $\phi \in \mathcal D(\R_+)$. Using then the change of variables $(X,V)=\Psi^{-t}(x,v)$ leads to 
\[
 \int_{\Omega\times\R^3}\chi(x,v)\left(\int_{\R_+}\left[f^\sharp \phi'(t) + h^\sharp \phi(t)\right]\ud t\right)\ud x\ud v=0
\]
and then to the announced result, since this is true for every $\chi \in \mathcal D(\Omega\times \R^3)$ and $\phi \in \mathcal D(\R_+)$.~~~$\square$

Equipped with this lemma, it is now possible to apply all the strategy developped in \cite{dolbeault}, and we obtain the following theorem:

\begin{thm}
Let $\Omega$ be either $\R^3$ or a regular subset of $\R^3$. Let the collision kernel be such that
\be
\label{hyp1}
0\le b\in L^1(\R^3\times S^2),
\ee 
and let 
\be
\label{hypf0}
f_0 \in L^\infty(\Omega\times\R^3), \qquad0\leq f_0 \leq 1.
\ee
Then, the problem \eqref{BFD}-\eqref{initial}-\eqref{BC} has a unique solution $f$ satisfying
\[
 f\in L^\infty(\R_+\times\Omega\times\R^3),\qquad 0\leq f \leq 1 ~~\textrm{a.e.}.
\]
Moreover, $f$ is absolutely continuous with respect to $t$.
\end{thm}
\begin{proof}[Sketch of proof]
 The main idea is to show that the function
\[
 T:f\mapsto f_0(\Psi^{-t}(x,v))+\int_0^t Q(\bar f)(s,\Psi^{s-t}(x,v))\ud s 
\]
has a fixed point, with
\[
 \bar f =\begin{cases}
	
           &0 ~~\textrm{if} ~~f\leq 0\\
	   &f ~~\textrm{if} ~~0\leq f\leq 1\\
	   &1 ~~\textrm{if} ~~f\geq 1
         \end{cases}
\]
This is achieved by showing that $T$ is contractive in $L^\infty([0,\theta]\times\Omega\times\R^3)$ if $\theta$ is small enough. Then, since
\[
 -B \max(f,0)\leq -B \bar f \leq Q(\bar f)\leq B(1-\bar f)\leq B(1-\min(1,f))
\]
where $B=\|b\|_{L^1(\R^3\times S^2)}$, it comes
\[
 -B \max(f^\sharp,0)\leq \partial_t f^\sharp \leq B(1-\min(1,f^\sharp)),
\]
which ensures that $f=\bar f$. It is then possible to reiterate the process on $[\theta,2\theta]$, and so on, to construct a global solution.
\end{proof}

\section{Conservation laws}

Under the (physically relevant) assumption on the collision kernel \eqref{hyp} the collision operator $Q(f)$ features the very interesting symmetry properties \eqref{symmetry}, thanks to which $f$ preserves some macroscopic quantities such as the total mass and the total kinetic energy. Indeed the following result holds as in \cite{dolbeault} or \cite{lu}:

\begin{thm}
\label{globalcons}
 Let $f_0$ satisfy \eqref{hypf0}, and $b$ satisfy \eqref{hyp} and \eqref{hyp1}. Then
\begin{itemize}
 \item if $f_0\in L^1(\Omega\times\R^3)$, then the solution $f$ to \eqref{BFD}-\eqref{initial}-\eqref{BC} belongs to $C^0(\R_+;L^1(\Omega\times\R^3))$ and
\[
 \iint_{\Omega\times\R^3}f(t,x,v)\ud x\ud v= \iint_{\Omega\times\R^3}f_0(x,v)\ud x\ud v,\qquad\forall t\in \R_+;
\]
 \item if $\iint_{\Omega\times\R^3}|v|^2 f_0\ud x\ud v <+\infty$, then the solution $f$ to \eqref{BFD}-\eqref{initial}-\eqref{BC} is such that the function $(t,x,v)\mapsto |v|^2 f(t,x,v)$ belongs to $C^0(\R_+;L^1(\Omega\times\R^3))$ and
\[
 \iint_{\Omega\times\R^3}|v|^2f(t,x,v)\ud x\ud v= \iint_{\Omega\times\R^3}|v|^2f_0(x,v)\ud x\ud v,\qquad\forall t\in \R_+.
\]
Moreover, the function $(t,x,v)\mapsto |v|^2 Q(f)(t,x,v)$ belongs to $L^\infty(\R_+;L^1(\Omega\times\R^3))$.
\end{itemize}

\end{thm}
\begin{proof}[Sketch of proof] Since the solution $f$ to \eqref{BFD} belongs to $C^0(\R_+;L^1(\Omega\times\R^3))$ it is easy to see that $Q(f) \in C^0(\R_+;L^1(\Omega\times\R^3))$ and the first point is a consequence of Fubini's theorem. The second point is more tricky and is dealt with as in \cite{dolbeault}, showing first the result with a fixed point in $L^\infty([0,\theta];L^1(\ud x (1+|v|^2)\ud v))$ assuming that $\int_{\R^3}\int_{S^2}b(w,\omega)|w|^2\ud \omega\ud w<\infty$ and then using a stability argument to relax this assumption.
\end{proof}

We can go further and show that even the microscopic quantities are conserved. This is a consequence of the following dispersion estimate coming from \cite{perthame}:

\begin{thm}
\label{thmperthame}
 Let $f\in L^1([0,T];L^1(\Omega\times\R^3))$ be solution of
\[
 \partial_t f + v.\nabla_x f =g,\qquad f(t=0)=f_0
\]
with boundary condition \eqref{BC}. Assume that
\[
 \iint_{\Omega\times\R^3}|v|^3 f_0(x,v)\ud x\ud v <+\infty\qquad \textrm{and}\qquad\int_0^T \iint _{\Omega\times \R^3}|v|^2 g(t,x,v)\ud x\ud v\ud t <C_0.
\]
Then, for any bounded subset $K$ of $\Omega$, we have
\[
 \int_0^T \ud t \int_{K}\ud x \int_{\R^3} |v|^3 f(t,x,v)\ud v \leq C_K.
\]
\end{thm}

This dispersion estimate is the key to obtain that $f$ locally conserves the mass, momentum and kinetic energy:

\begin{cor}
\label{theorem}
 Let the collision kernel $b$ satisfy \eqref{hyp} and \eqref{hyp1}. Assume that the initial data satisfy \eqref{hypf0}, and
\[
 \iint_{\Omega\times\R^3}(1+|v|^3)f_0(x,v)\ud x\ud v <+\infty.
\]
Then, the solution $f$ to \eqref{BFD}-\eqref{initial}-\eqref{BC} satisfies, in distributional sense, the local conservation laws \eqref{conslaw}.
\end{cor}
Note that because of the boundary condition \eqref{boundary}, it is easy to see that the momentum is tangential to the boundary:
\[
 n.\int_{\R^3}vf\ud v=0\qquad\quad\forall x\in\partial\Omega.
\]
\textbf{Proof.}
In the sense of distributions, $f$ satisfies \eqref{BFD}. Let $\phi \in \mathcal D(\R_+\times\Omega)$ and $\Psi_R\in \mathcal D(\R^3)$ satisfy:
\[
 \Psi_R(v)=
\begin{cases}
  &1\quad\textrm{if}~|v|\leq R\\
	 &0\quad \textrm{if}~|v|>2R.
 \end{cases}
\]
Then equation \eqref{BFD} implies
\[
\begin{split}
 &\int_{\R_+\times\Omega}\partial_t\phi \left[\int_{\R^3}\Psi_R (v)\begin{pmatrix} 1 \\v\\|v|^2\end{pmatrix}f\ud v \right] \ud x\ud t+ \int_{\R_+\times\Omega}\nabla_x \phi\left[\int_{\R^3}\Psi_R(v) \begin{pmatrix}v\\v\otimes v\\v|v|^2\end{pmatrix} f\ud v\right]\ud x\ud t\\
&=- \int_{\R_+\times\Omega}\phi \left[\int_{\R^3}\Psi_R(v) \begin{pmatrix} 1\\ v\\|v|^2\end{pmatrix} Q(f)\ud v\right]\ud x\ud t
\end{split}
\]
We may pass to the limit $R\to +\infty$ thanks to theorems \ref{globalcons} and \ref{thmperthame}. At the end, we get the desired result.~~$\Box$

\end{document}